\documentclass[12pt]{amsart} 

\usepackage{amssymb} 
\usepackage{amsmath} 
\usepackage{fullpage} 
\usepackage{url}

\newtheorem{theo}{Theorem}

\theoremstyle{definition}
\newtheorem{defi}{Definition}
\theoremstyle{remark}
\newtheorem{rem}{Remark}

\def\Er{{\mathbb E}}

\def\Pr{{\mathbb P}}

\def\Rr{{\mathbb R}}

\def\Bc{{\mathcal{B}}}

\def\Fc{{\mathcal{F}}}

\def\one{{\rm \bf 1}}
\def\({\left(}     
\def\){\right)}    
\def\[{\left[}     
\def\]{\right]}

\def\as{{\frenchspacing a.s.}~}
\begin{document}
\title{Commonotonicity and $L^1$ Random Variables} 
\author{Freddy Delbaen}
\address{Departement f\"ur Mathematik, ETH Z\"urich, R\"{a}mistrasse
101, 8092 Z\"{u}rich, Switzerland}
\address{Institut f\"ur Mathematik,
Universit\"at Z\"urich, Winterthurerstrasse 190,
8057 Z\"urich, Switzerland}
\date{First version May 1, 2019, this version \today}

\begin{abstract}
It is proved that in suitable filtrations every pair of integrable random variables is the conditional expectation of a pair of commonotone integrable random variables. \end{abstract}

\maketitle

\section{Notation}

We will work with a probability space equipped with three sigma algebras $(\Omega,\Fc_0\subset \Fc_1\subset \Fc_2,\Pr)$.  The sigma algebra $\Fc_0$ is supposed to be trivial $\Fc_0=\{\emptyset,\Omega\}$ whereas the sigma algebra $\Fc_2$ is supposed to express innovations with respect to $\Fc_1$.  Since we do not put topological properties on the set $\Omega$ we will make precise definitions later that do not use conditional probability kernels.  But essentially we could say that we suppose that conditionally on $\Fc_1$ the probability $\Pr$ is atomless on $\Fc_2$. In a previous paper we give equivalent conditions for this property, \cite{FDcom}.   The space $L^\infty(\Fc_i)$ will denote the space of bounded $\Fc_i$ measurable random variables, modulo almost sure equality \as.   The space $L^1(\Fc_i)$ is the space of integrable $\Fc_i$-measurable random variables, modulo equality \as.

We say that two random variables $\xi,\eta$ are commonotone if there are two nondecreasing functions $f,g\colon \Rr\rightarrow \Rr$ and a random variable $\zeta$ such that $\xi=f(\zeta), \eta=g(\zeta)$.  {\it Commonotonicity can be seen as the opposite of diversification}. Commonotone pairs play an important role in risk analysis, e.g. actuarial mathematics, finance, game theory, \ldots . If $\zeta$  increases then both $\xi$ and $\eta$ increase (or better do not decrease). By the way in case $\xi$ and $\eta$ are commonotone then one can choose $\zeta=\xi+\eta$, see \cite{FDbook}.  It can be shown (an exercise) that in this case one can choose representatives (still denoted $\xi,\eta$ such that $(\xi(\omega)-\xi(\omega'))(\eta(\omega)-\eta(\omega'))\ge 0$ for all $\omega,\omega'$.  

We say that a set $E\subset \Rr^2$ is commonotone  if $(x,y),(x',y')\in E$ implies $(x-x')(y-y')\ge 0$.  Using this, random variables $\xi,\eta$ are commonotone if and only if, the support of the image measure of $(\xi,\eta)$ is a commonotone set.  If $E$ is a commonotone set then also the closure $\overline{E}$ is  commonotone.  To get commonotone sets we can use the following technique.  Given a subset $I\subset \Rr$ and two nondecreasing functions $g,h\colon I\rightarrow \Rr$, we can put $E=\{(g(t),h(t)\mid t\in I\}$.  The reader can make pictures when for instance $I$ is an interval or $I=\Rr$.  In this case $E$ becomes a nondecreasing curve in $\Rr^2$, a typical example of a commonotone set.

\section{A Special Commonotone Set in $\Rr^2$}

To deal with integrable random variables we use the same technique as in \cite{FDcom}.  The unboundedness of the random variables poses some technical problems that we will overcome by using a special commonotone set in $\Rr^2$.  The construction seems a little bit complicated but the reader can make a drawing to see what happens.  The set will be constructed by induction.   The first step consists in taking the curve obtained as the concatenation of the convex intervals that join the points $$(-4,-4)\rightarrow (-4,-2)\rightarrow (0,0)\rightarrow (4,2)\rightarrow (4.4).$$
The convex hull of this set is a parallelogram $P_1$, with parallel vertical sides given by the segments
$$((-4,-4)\rightarrow(-4,-2)\text{  and   } (4,2)\rightarrow (4,4).$$
Note that every point of $P_1$ is the convex combination of points taken on the vertical sides.  An easy and continuous way to obtain such convex combination goes as follows.  Through a point in $P_1$ take a line parallel to the ``skew" sides of $P_1$ and see where it intersects the vertical sides.  We will also need some norm estimates.  Since all norms on $\Rr^2$ are equivalent, we will use the one that avoids constants $C$ that {\it change from one line to the next or even within the same line}.  The set $P_1$ is symmetric around the origin, is convex and compact and  can be seen as the unit ball of a norm. Hence $\Vert (x,y)\Vert \le 1$ if and only if $(x,y)\in P_1$.

Before we describe the general step, let us see what happens on the second step.  For expository reasons we extend the curve made in step 1 by adding two more pieces.  One at the positive side and one at the negative side:
$$
(-8,-8)\rightarrow (-8,-4)\rightarrow (-4,-4)
$$ 
and
$$
(4,4)\rightarrow (8,4)\rightarrow (8,8).
$$
The parallelogram $P_2$ is the convex hull of the line segments
$$(-8,-8)\rightarrow (-8,-4)\text{  and   }(8,4)\rightarrow (8,8)$$
In doing so $P_2$ becomes the double of $P_1$, i.e. $P_2=2 P_1=\{(x,y)\mid \Vert(x,y)\Vert\le 2\}$.   Therefore $P_1\subset P_2$. If $(x,y)\in P_2\diagdown P_1$ we have $\Vert (x,y)\Vert \ge 1$.

The general step is now clear.  At stage $n-1$ we have a curve from $(-2^n,-2^n)$ to $(2^n,2^n)$, a parallelogram $P_{n-1}$ and we can construct extra lines
$$
(-2^{n+1},-2^{n+1})\rightarrow (-2^{n+1},-2^n)\rightarrow (-2^n,-2^n)
$$ and 
$$
(2^n,2^n)\rightarrow (2^{n+1},2^n)\rightarrow (2^{n+1},2^{n+1})
$$
as well as a parallelogram $P_n$ with $P_n=2 P_{n-1}=2^{n-1}P_1$, $P_{n-1}\subset P_n$. $P_n=\{(x,y)\mid\Vert(x,y)\Vert\le 2^{n-1}\}$.  If $(x,y)\in P_n\diagdown P_{n-1}$ we have $\Vert (x,y)\Vert\ge 2^{n-2}$.  Using the same procedure as for $P_1$ we can represent each point in $(x,y)\in P_n\diagdown P_{n-1}$ as a convex combination of two points $(u_1,v_1);(u_2,v_2)$ on the vertical sides of $P_n$.  Important for later use are the inequalities $\Vert (u_i,v_i)\Vert\le 2\Vert (x,y)\Vert$. The union of all the curves used to construct the different parallelograms is denoted by $E$.  $E$ is a commonotone set.

Pasting together all these domains, $P_n\diagdown P_{n-1}, n\ge 2$, $P_1$ and the convex combinations defined on them, gives us Borel measurable functions
\begin{enumerate}
\item $\lambda\colon \Rr^2\rightarrow [0,1]$
\item $(u_1,v_1)\colon \Rr^2\rightarrow E$, $(u_2,v_2)\colon \Rr^2\rightarrow E$
\item For $\Vert (x,y)\Vert > 1$ we have $\Vert (u_i,v_i)\Vert\le 2 \Vert (x,y)\Vert$
\item For $\Vert (x,y)\Vert\le 1$ we have $\Vert (u_i,v_i)\Vert  = 1$.
\item for all $(x,y)\in \Rr^2$: $(x,y)=\lambda (u_1,v_1)+(1-\lambda)(u_2,v_2)$
\end{enumerate}

\section{Atomless Extension}                                                                                                                                                                                                                                                                  

\begin{defi}  We say that $\Fc_2$ is atomless conditionally to $\Fc_1$ if the following holds.  If $A\in \Fc_2$ then there exists a set $B\subset A$, $B\in\Fc_2$, such that $0< \Er[\one_B\mid\Fc_1]<\Er[\one_A\mid\Fc_1]$ on the set $\{\Er[\one_A\mid\Fc_1]>0\}$.
\end{defi}
In  case the conditional expectation could be calculated with a -- under extra topological conditions -- regular probability kernel, say $K(\omega, A)$, then the above definition is a measure theoretic way of saying that the probability measure $K(\omega, .)$ is atomless for almost every $\omega\in\Omega$.
In \cite{FDcom} the following theorem is proved.
\begin{theo} Are equivalent:
\begin{enumerate}
\item $\Fc_2$ is atomless conditionally to $\Fc_1$
\item For $A\in\Fc_2$ there is $B\subset A$ such that $\Pr\[0< \Er[\one_B\mid\Fc_1]<\Er[\one_A\mid\Fc_1]\]>0$.
\item  There exists an atomless sigma-algebra $\Bc\subset \Fc_2$ that is independent of $\Fc_1$.
\item There is an increasing family of sets $(B_t)_{t\in [0,1]}$ such that $\Er[\one_{B_t}\mid\Fc_1]=t$.  The sigma algebra $\Bc$, generated by the family $(B_t)_t$ is independent of $\Fc_1$.  The system $(B_t)_t$ can also be described as $B_t=\{U\le t\}$ where $U$ is a random variable that is uniformly distributed on $[0,1]$,  $U$ and $\Fc_1$ are independent.
\item There is a uniformly $[0,1]$ distributed random variable $U\colon \Omega\rightarrow [0,1]$, independent of $\Fc_1$, such that for every $h\colon\Omega \rightarrow [0,1]$ which is $\Fc_1$ measurable we have $\Er[\one_{\{U\le h\}}\mid \Fc_1]=h$
\end{enumerate} 
\end{theo}

\section{The main Result}
This section is devoted to the proof of
\begin{theo}  Suppose  that $\Fc_2$ is atomless conditionally to $\Fc_1$. For any two integrable $\Fc_1$ measurable random variables, $f,g$, we can find two commonotone $\Fc_2$ random variables $\xi,\eta$ such that $f=\Er[\xi\mid\Fc_1], g=\Er[\eta\mid\Fc_1]$.  The (2-dimensional) random variable $(\xi,\eta)$ has the same tail behaviour as $(f,g)$ and is therefore integrable.  More precisely for the norm with unit ball $P_1$ we have almost surely $\Vert (\xi,\eta)\Vert\le \max(2\Vert (f,g)\Vert, 1)$.
\end{theo}
{\bf Proof } For given integrable $\Fc_1$ measurable $(f,g)\colon\Omega\rightarrow \Rr^2$, we define  (using the notation of the end of section 2)  $\Lambda = \lambda\circ (f,g)$, $(U_1,V_1)=(u_1,v_1)\circ (f,g)$, $(U_2,V_2)=(u_2,v_2)\circ (f,g)$.  Clearly these functions are $\Fc_1$ measurable.  We now define $(\xi,\eta)=\one_{\{U\le\Lambda\}}(U_1,V_1)+\one_{\{U>\Lambda\}}(U_2,V_2)$, where $U$ is defined as in the previous section.  Taking conditional expectations gives the desired result.  The inequalities follow from the corresponding inequalities for $(u_1,v_2)$ and $(u_2,v_2)$.
\begin{rem} Conditional expectations with respect to $\Fc_1$ can be defined using localisation to sets in $\Fc_1$. This allows to extend the definition to cases where $(f,g)$ are not necessarily integrable.  When needed a reader can do this according to the need.  Adding a theorem that covers {\it all } cases would lead to an even more  psychedelic formulation -- something we want to avoid.
\end{rem}
\begin{rem}  We can say a little bit more than just integrability.  If $(f,g)$ is in a solid space then $(\xi,\eta)$ is in the same space.  This can be applied to $L^p$ spaces, Orlicz spaces, Orlicz hearts and ordered spaces coming from monetary utility functions. We do not give applications or details since they fall outside the scope of this short note.
\end{rem}

\end{document}